\documentclass[12pt]{article}
\usepackage{amsmath}
\usepackage{amscd}
\usepackage{amssymb}
\begin{document}

\newcommand{\R}{\mathbb{R}}
\newcommand{\C}{\mathbb{C}}
\newcommand{\g}{\mathcal{G}}
\newcommand{\A}{{\mathcal{A}}}
\newcommand{\I}{{\mathcal{I}}}
\newcommand{\lc}{{\mathcal{L}}}
\vskip 2cm
  \centerline{\Large \bf Deformation Quantization of Coadjoint Orbits}

  \bigskip

\centerline{ M. A. Lled\'o}
\bigskip

\centerline{\it  Dipartimento di Fisica, Politecnico di Torino,}
\centerline{\it Corso Duca degli Abruzzi 24, I-10129 Torino, Italy, and}
\centerline{\it INFN, Sezione di Torino, Italy.}
\centerline{{\footnotesize e-mail: lledo@athena.polito.it}}

\bigskip

\begin{abstract}
A method for the deformation quantization of coadjoint orbits
of semisimple Lie groups is proposed.  It is based on the algebraic
structure of the orbit. Its relation to geometric
quantization and differentiable deformations is explored.
\end{abstract}

Let $G$ be  a complex Lie group of dimension $n$ and $G_R$ a real form of $G$.
Let $\g$ and $\g_R$ be their respective Lie algebras with Lie bracket $[\;,\;]$.
 As it is well known, $\g_R^*$ has a Poisson structure,
\begin{equation}
\{f_1,f_2\}(\lambda)=<[(df_1)_{\lambda},(df_2)_{\lambda}],\lambda>,
\qquad f_1,f_2 \in C^{\infty}(\g_R^*), \quad \lambda \in \g_R^*.
\label{pb}
\end{equation}
Choosing a basis $\{X_1,\dots X_n\}$ of $\g_R$ and its dual, $\{\xi^1,\dots
\xi^n\}$, the Poisson bracket can be written as
$$
\{f_1,f_2\}(x)=c_{ij}^k\lambda_k\frac{\partial f_1}{\partial x_i}
\frac{\partial f_2}{\partial x_j},\qquad x=\sum_{i=1}^nx_i\xi^i\in \g^*.
$$
Notice that this Poisson bracket is never symplectic, in particular it
is 0 at the origin.
Under the action of $g\in G_R$, it satisfies
$g^*\{f_1,f_2\}=\{g^*f_1,g^*f_2\}$, so $G_R$ is a group of automorphisms
of the Poisson algebra $C^\infty(\g_R^*)$. The action
of $G_R$ on $\g_R^*$ is not transitive, so $\g_R^*$ is foliated in orbits.
This foliation coincides with the foliation given by the Hamiltonian
vector fields of (\ref{pb}). So the orbits of the coadjoint action of a
Lie group are symplectic manifolds.

We want to describe formal deformations of the Poisson
algebra $C^\infty(\Theta)$ ($\Theta$ is a coadjoint orbit) or of some
subalgebra of it. It is convenient
to work with the complexification of the Poisson
algebra.

An associative algebra $\A_h$ over $\C[[h]]$ is a
formal deformation of
a Poisson  algebra $(\A, \{\;,\;\})$ over $\C$ if there exists an isomorphism of  
$\C[[h]]$-modules $\psi: \A[[h]]\longrightarrow \A_h $ satisfying the
following properties:

\smallskip

\noindent {\bf a.} $\psi^{-1}(F_1F_2)=f_1f_2$ mod($h$) where $F_i\in \A_h$ are
such that $\psi^{-1}(F_i)=f_i$ mod$(h)$, $f_i\in \A$. (By mod$(h)$ we mean that the
projections $p:\A[[h]]\longrightarrow \A[[h]]/h\A[[h]]$ of both
quantities coincide).

\smallskip

\noindent {\bf b.} $\psi^{-1}(F_1 F_2 - F_2 F_1) = h\{f_1,f_2\}$ mod($h^2$).

\smallskip

An example of interest for our purposes is the polynomial algebra on
$\g^*$.  A formal deformation of $\mbox{Pol}(\g^*)$ \cite{ho} is given by the
algebra $U_h= T_{\C[[h]]}(\g)/\lc_h$, where $T_{\C[[h]]}(\g)$
is the tensor algebra over $\C[[h]]$ and $\lc_h$ is the
  proper two sided ideal 
$$
{\lc}_h=\sum_{X,Y \in \g}
T_{\C[[h]]}(\g) \otimes(X \otimes Y - Y \otimes X - h[X,Y])
\otimes T_{\C[[h]]}(\g)\subset T_{\C[[h]]}(\g).
$$
The isomorphism $\psi:\mbox{Pol}(\g^*)\longrightarrow U_h$ is not
canonical. A possible choice is in terms of a Poincar\'e-Birkhoff-Witt basis,
\begin{equation}
\psi(x_{i_1}x_{i_2}\cdots x_{i_k})=X_{i_1}\cdot X_{i_2}\cdots X_{i_k},
\qquad 1\leq i_1\leq\cdots \leq i_k\leq n.
\label{sp1}
\end{equation}
Another choice is  the symmetrizer map,
\begin{equation}
\hbox{Sym}(x_{i_1}x_{i_2}\cdots x_{i_k})=\frac{1}{k!}\sum_{\sigma \in S_k}X_{\sigma(i_1)}
\cdot X_{\sigma(i_2)}\cdots X_{\sigma(i_k)},
\label{sp2}
\end{equation}
where $S_k$ is the group of permutations of order $k$.

Given a choice for $\psi$ one can define an associative product (star
product) on $\A[[h]]$ by
$$
a\star_\psi b=\psi^{-1}(\psi(a)\cdot \psi(b)).
$$
Then, for any choice of $\psi$, $(\A_h,\star_\psi)$ is an algebra isomorphic
to $\A_h$. With the star product we recover the semiclassical
interpretation of the elements of the algebra as functions on the phase space.
 The star product can always be written
as a formal series
$$a\star_\psi b=ab +\sum_{n>0}h^nC_\psi^n(a,b),$$
where $C_\psi^n$ are some bilinear operators. Let  $\star$ and 
$\star'$ be two isomorphic star products, 
$$
a\star b= T^{-1}(T(a)\star T(b)), \qquad T:\A[[h]]\longrightarrow A[[h]].
$$
It is clear that $T$ can be written as
$$
T(a)=\sum_{n\geq 0} h^nT^n(a),
$$
and because of property {\bf a}, $T^0$ must be an automorphism of the
commutative algebra $\A[[h]]$. If $T_0$ is the identity we say that
$\star$ and $\star'$ are  equivalent (or gauge equivalent) star
products. (\ref{sp1}) and (\ref{sp2}) are two equivalent star products.

For $\A$ being the full algebra of $C^\infty$ functions on the Poisson manifold,
if the operators $C_\psi^n$ are bidifferential operators we say that the star
product is differentiable. Gauge equivalence can be restricted to the
class of differentiable star products by considering only differentiable
$T^n$. Notice that the differentiability is a property of the particular
star product and not of the formal deformation. The star products
(\ref{sp1}) and (\ref{sp2}) can be extended to $C^\infty(\Theta)$ as differentiable
star products, but we will see later an
example of a star product corresponding to the same formal deformation
which is not differentiable \cite{cg}.

\medskip

We will consider only semisimple Lie groups. The semisimple coadjoint orbits
of $G$ on $\g$ are complex algebraic varieties defined over $\R$. They
are given by the invariant polynomials. If $l$ is the rank of $G$, we
can choose $l$ homogeneous polynomials 
$p_i(x)$, $i=1,\dots , l$ generating the subalgebra of invariant
polynomials, $\C[p_1,\dots p_n]$. Then the semisimple coadjoint orbits are given by
the algebraic equations
\begin{equation}
p_i(x)=c_i.
\label{av}
\end{equation}
(see for example \cite{va1}). The intersection of the complex orbit with
$\g_R$ is a real algebraic variety  consisting on a finite number of
connected components, which are
orbits of the real form of the group. For the  compact real form
there is only one connected component.

It is easy to check that the star products  \ref{sp1} and \ref{sp2} do not
restrict well to the orbit, that is, in general
$$
a\star p_i|_\Theta\neq 0.
$$
We want to know if there is some choice of $\psi$ that gives a star product
which restricts to $\Theta$.

In the approach of geometric quantization, the algebra of quantum observables
is given by the quotient of $U_h$ by a certain ideal. This ideal, $I_h$, is prime
and $\mbox{Ad}_G$-invariant, so there is a well defined action of $G$ on
$U_h/I_h$.  We summarize here the results of \cite{fl}, where the
quantization of the coadjoint orbits are obtained in terms of the
quotient of the enveloping algebra by a prime, Ad$_G$-invariant ideal.
We consider the polynomial algebra over the real algebraic manifold
(union of orbits) defined by (\ref{av}),
$$
\mbox{Pol}(\Theta)=\mbox{Pol}(\g^*)/\I_0,\qquad \I_0=\{p\in \mbox{Pol}(\g^*)/p|_\Theta=0\}.
$$
This is the Poisson algebra that we want to deform. We quote first a
result from Varadarajan \cite{va2} that we need.

\smallskip

{\bf Lemma (1)}. {\it Let $x\in\g^*$ be a  regular element of $\g^*$
(or equivalently, a
point in which the centralizer has dimension equal to the rank of $\g^*$).
Then $(dp_1)_{x}$, ..., $(dp_l)_{x}$ are linearly independent.}

\smallskip

From now on we will restrict to regular orbits only. It is clear that in
this case $\I_0$ is generated by $p_1-c_1,\dots p_l-c_l$. We consider now
the elements in $U_h$ that are the image of $p_i$ by the symmetrizer,
$P_i=\mbox{Sym}(p_i)$, called Casimir operators.
 $P_i$ are central elements in $U_h$ and they are
also Ad$_G$-invariant (the symmetrizer
commutes with the action of $G$).

Let $\I_h$ be the ideal generated by
$P_i-C_i(h)$, where $C_i(0)=c_i$. Then $U_h/\I_h$ is a formal deformation
of $\mbox{Pol}(\Theta)$. The technical assumption of regularity is needed
to prove the existence of a $\C[[h]]$-module isomorphism
$\psi:\mbox{Pol}(\Theta)[[h]]\longrightarrow U_h/\I_h$, which is not obvious.
The ideal $\I_h$ itself is Ad$_G$-invariant, so $G$ has a natural action
by automorphisms on the algebra $U_h/\I_h$. For special values of
$C_i(h)$, $\I_h$ is in the kernel of an irreducible unitary
representation of $G_R$.

This deformation of polynomials can be
specialized for any value of $h$. For SU(2),
$$
[ H, X]=\hbar 2 X,\quad [ H, Y]=-\hbar 2Y,
 \quad [ X, Y]=\hbar  H,
$$
the Casimir operator is
$$
 P=\frac{1}{2}(X Y+ Y X +\frac{1}{2} H^2).
$$
It was  shown in \cite{fl} that with the choice
$$
C=l(l+\hbar )),\qquad l=\hbar {m/2},
$$
the algebra obtained is the same than the one obtained in geometric quantization.

\medskip
According to our definition, a star product in Pol$(\Theta)$ is given by
a $\C[[h]]$-module isomorphism  $\tilde\psi:\mbox{Pol}(\Theta)[[h]]\longrightarrow
U_h/\I_h$. In particular, to obtain a star product in Pol$(\g^*)$ which
restricts to the orbit one should look for an isomorphism
$\psi:\mbox{Pol}(\g^*)[[h]]\longrightarrow U_h$ such that the following
diagram commutes
\begin{equation}
\begin{CD}
\mbox{Pol}(\g^*)[[h]]@>\psi>>U_h\\
@VV{\pi}V      @VV{\hat \pi}V\\
\mbox{Pol}(\Theta)[[h]]@>\tilde\psi>>U_h/I_h.
\end{CD}
\label{dia}
\end{equation}
$\pi$ and $\hat\pi$ are the natural projections. If $\psi(\I_0)=\I_h$,
then $\tilde \psi$ is defined uniquely by the diagram (\ref{dia}). An
example of such star product is given in \cite{cg}, where it is
also shown that it is not differentiable. Since $\hat\psi$ is not unique
one could ask if there is some choice that renders it differentiable.
This issue will be addressed in \cite{fll}.


\begin{thebibliography}{99}

\bibitem{ho} J. Hoppe, {\it Quantum theory of a massless relativistic
 surface and a two dimensional bound state problem}.
MIT PhD thesis, (1982).

\bibitem{cg} M. Cahen and S. Gutt, {\it C.R. Acad. Sc. Paris}  {\bf 296} (1983) 821-823;
in {\sl Gravitation
and Cosmology}, Monographs and Textbooks in Physical Sciences,  Bibliopolis.
 Eds W. Rundler and A. Trautman, (1987);
 {\it Bull. Soc. Math. Belg.} {\bf 36 B}, 207-221, (1987).


\bibitem{va1} V. S. Varadarajan,  {\it Harmonic  Analysis on Real Reductive
Groups.} Lecture Notes in Mathematics, no. 576. Springer-Verlag, (1977).


\bibitem{fl} R. Fioresi and M. A. Lled\'o, math.QA/9906104. To appear in
the {\it Pacific Journal of Mathematics}.

\bibitem{va2} V. S. Varadarajan,
{\it Amer. J. Math}, {\bf 90}, (1968).


\bibitem{fll} R. Fioresi, A. Levrero and M. A. Lled\'o. In preparation.

\end{thebibliography}
\end{document}